\documentclass[leqno,11pt]{article}
\usepackage{latexsym}
\usepackage{amssymb}
\usepackage{amsfonts}
\usepackage{amsmath}
\usepackage{color}
\usepackage{pifont}

\definecolor{verde}{rgb}{0.30,0.7,0.00}
\definecolor{bleu}{rgb}{0.00,0.4,0.90}

\makeatletter
\long\def\unmarkedfootnote#1{{\long\def\@makefntext##1{##1}\footnotetext{#1}}}
\makeatother

\newtheorem{definition}{Definition}[section]

\newtheorem{theorem}[definition]{Theorem}

\newtheorem{corollary}[definition]{Corollary}

\def\R{\mathbb R}

\def\rn{{{\R}^n}}
\def\rN{{{\R}^N}}

\def\r2{{{\R}^2}}

\def\O{\Omega}

\def\to{\rightarrow}

\def\bdm#1\edm{\begin{displaymath}#1\end{displaymath}}
\def\be#1\ee{\begin{equation}#1\end{equation}}
\def\barr#1\earr{\begin{align}#1\end{align}}

\newcommand{\medint}{-\kern  -,395cm\int}
\newcommand{\medintinrigo}{-\kern  -,315cm\int}
\newcommand{\medelle}{-\kern  -,235cm L}
\newcommand{\medellenrigo}{-\kern  -,180cm L}
\newcommand{\qed}{\thinspace\null\nobreak\hfill
\hbox{\vbox{\kern-.2pt\hrule height.2pt
depth.2pt\kern-.2pt\kern-.2pt \hbox to1.8mm {\kern-.2pt\vrule
width.4pt \kern-.2pt\raise1.8mm\vbox to.2pt{} \lower0pt\vtop
to.2pt{}\hfil\kern-.2pt \vrule
width.4pt\kern-.2pt}\kern-.2pt\kern-.2pt \hrule height.2pt depth.2pt
\kern-.2pt}}\par\medbreak}

\title{A symmetrization result for a class of\\
 anisotropic elliptic problems}
\author{A. Alberico\thanks{Istituto per le Applicazioni del Calcolo ``M.
Picone'' (I.A.C), Sez. Napoli, Consiglio Nazionale delle Ricerche
(C.N.R.), Via P. Castellino 111, 80131 Napoli, Italy.
E--mail:a.alberico@na.iac.cnr.it} -- G. di
Blasio\thanks{Dipartimento di Matematica e Fisica,Universit\`{a}
degli Studi della Campania \textquotedblleft Luigi Vanvitelli
\textquotedblright, Via Vivaldi, 43 - 81100 Caserta, Italy. E--mail:
giuseppina.diblasio@unicampania.it} -- F. Feo\thanks{Dipartimento di
Ingegneria, Universit\`{a} degli Studi di Napoli \textquotedblleft
Pathenope\textquotedblright, Centro Direzionale Isola C4 80143
Napoli, Italy. E--mail: filomena.feo@uniparthenope.it}}
\date{}

\begin{document}

\maketitle

\begin{abstract}
%In this paper we prove some estimates  for weak  solutions to
%anisotropic elliptic problems, whose prototype is
%\[
%\left\{
%\begin{array}
%[c]{ll}%
%-\displaystyle\sum_{i=1}^N \frac{\partial}{\partial x_{i}}\left(
%\left\vert \frac{\partial u}{\partial x_{i}}\right\vert
%^{p_{i}-2}\frac{\partial u}{\partial
%x_{i}}\right)+b(u)=f(x) & \hbox{in $\Omega$}\\
%& \\
%u=0 & \hbox{on $\partial\Omega$},
%\end{array}
%\right.
%\]
%where $\Omega$  is a bounded open subset of $\R^{N}$, with Lipschitz
%continuous boundary, $N\geq2$, $1\leq p_{i}<\infty$ for
%$i=1,\ldots,N$ such that their harmonic mean is greater than $1$,
%$b$ is a continuous, non-decreasing function such that $b(0)=0$ and
%$f$ is a nonnegative function with a suitable summability.
We prove estimates for weak solutions to a class of Dirichlet
problems associated to
anisotropic elliptic equations with a zero order term..% As a consequence, a priori
%estimates for norms of the relevant solutions are derived.
\end{abstract}

\bigskip


\numberwithin{equation}{section}

\section{\textbf{Introduction}}

We consider the class of Dirichlet problems for anisotropic elliptic
equations, whose prototype has the form
\begin{equation}
\left\{
\begin{array}
[c]{ll}%
-\displaystyle\sum_{i=1}^N \left( \left\vert u_{x_i}\right\vert
^{p_{i}-2} u_{x_i}\right)_{x_i}+b\left(  u\right)  =f(x) & \hbox{in $\Omega$}\\
& \\
u=0 & \hbox{on $\partial\Omega$},
\end{array}
\right.  \label{proto2}%
\end{equation}
where $\Omega$ is a bounded open subset of $\R^{N}$ with Lipschitz
continuous boundary, $N\geq2$, $p_{i}\geq 1$ for $i=1,\ldots,N$ such
that their harmonic mean  $\overline{p}$ is greater than $1$, the
subscript $x_i$ denotes partial derivative with respect to $x_i$,
$b$ is a continuous, non-decreasing function such that $b(0)=0$ and
$f$ is
a nonnegative function with a suitable summability. %Here, $u_{x_i}$
%denotes the partial derivative of $u$ with respect to $x_i$.\

The anisotropy of problem (\ref{proto2})  depends on differential
operator whose growth with respect to the partial derivatives of $u$
is governed by different powers. In the last years anisotropic
problems have been extensively studied by many authors (see
\emph{e.g.} \cite{AdBF2,AdBF3,antontsev-chipot-08, BMS,
DiNardo-Feo-Guibe, DiNardo-Feo, FGK, FGL, FS, Gi, Mar}).
\par
The
growing interest has led to an extensive investigation also for
problems governed by fully
anisotropic growth conditions (see \cite{AC,A,AdBF1, Clocal,cianchi anisotropo}%
) and problems related to different type of anisotropy (see
\emph{e.g.} \cite{AFTL, BFK, DPdB,DPG}).

Our goal is to obtain  an estimate of concentration of a weak
solution to problem (\ref{proto2})  via symmetrization methods. The
use of the standard isoperimetric inequality in the study of
isotropic elliptic Dirichlet problems was introduced in \cite{Ma1,
Ma2} and independently  in \cite{Ta1, Ta2}. Variants and extensions
from these papers have been developed in a rich literature. We refer
to Vazquez \cite{V} and Trombetti \cite{T} for a quite comprehensive
bibliography on this and related topics.

 It is well known that when isotropic elliptic Dirichlet problems with a zero
 order term are considered,
 the situation is quite different
 if we assume or not a sign condition (see, \emph{e.g.}, \cite{Dz85, Dz, Md, v1, V}).
In the anisotropic setting there are two different cases as well.
Indeed, when $b(u)u\geq 0$, it is showed  (see, \emph{e.g.},
\cite{cianchi anisotropo}) that  the symmetric rearrangement of a
solution  $u$ to anisotropic problem (\ref{proto2}) is pointwise
dominated by the radial solution to an isotropic problem, defined in
a ball, with a radially symmetric decreasing data and with no zero
order term.
%If assumption $b(u)u\geq 0$ holds,
% rearrangement methods allow to obtain
%a pointwise comparison result between the solution $u$ to
%anisotropic problem (\ref{proto2}) and the radial solution to an
%isotropic simpler problem defined in a ball with a radially
%symmetric decreasing data with no zero order term (see \cite{cianchi
%anisotropo}).
Otherwise, with no sign condition on $b(u)u$, we prove  an integral
comparison result between a solution $u$ to anisotropic problem
\eqref{proto2} and the radial solution to a suitable isotropic
problem defined in a ball, with a radially symmetric decreasing data
again but, this time, which preserves a zero order term.

\par
 Just to give  an idea of our results, let us consider problem
(\ref{proto2}) when the domain $\Omega$\ is $B_{R}(0)$, the ball
centered at the origin and with radius $R>0$. We take into account
two smooth  strictly increasing functions $b$ and $\widetilde{b}$
having the same domain such that $b(0)=\widetilde{b}(0)=0$, and two
positive decreasing radial symmetric functions $f$ and
$\widetilde{f}$ defined in  $B_{R}(0).$ Denote by $b^{-1}$ and
$(\,\widetilde{b}\,)^{-1}$ the inverse function of $b$ and
$\widetilde{b}$, respectively. Suppose that
\[
(  (\,\widetilde{b}\,)^{-1})  ^{\prime}(s)\leq\left( b^{-1}\right)
^{\prime}(s)\qquad \hbox{for every}\;\;s\in \R
\]
and that the datum $f$ is less concentrated than the
datum $\widetilde{f},$ \textit{i.e.}%
\[
\int_{B_{r}(0)}f(x)\;dx\leq\int_{B_{r}(0)}\widetilde{f}(x)\;dx\qquad
\hbox{for every}\;\; 0\leq r\leq R.
\]
 Then, we are going to prove that%
\[
\int_{B_{r}(0)}b (u^{\bigstar}(x))\;dx\leq\int_{B_{r}(0)}\widetilde
{b} \left(\widetilde{u}(x)\right)\;dx\qquad \hbox{for
every}\;\;0\leq r\leq R,
\]
where $u^{\bigstar}$ is the symmetric decreasing rearrangement of
the solution $u$ to problem (\ref{proto2}) and $\widetilde{u}$ is
the solution to the
following problem%
\[
\left\{
\begin{array}
[c]{ll}%
-\operatorname{div}\left(  \left\vert \nabla\widetilde{u}\right\vert
^{\overline{p}-2}\nabla\widetilde{u}\right) +\widetilde{b}
(\widetilde{u})=\widetilde{f}(x) & \hbox{in $B_R(0)$}\\
& \\
\widetilde{u}=0 & \hbox{on $\partial B_R(0)$}.
\end{array}
\right.
\]
The paper is organized as follows. In Section 2 we recall some
backgrounds on the anisotropic spaces and on the properties of
symmetrization. In Section 3 we state our main results, proved in
Section 4.

\section{Preliminaries}

%\subsection{Anisotropic spaces}

Let $\Omega$ be a bounded open subset of $\R^{N}$, $N\geq2$, and let
$1\leq p_{1},\ldots,p_{N}<\infty$ be $N$ real numbers. The
anisotropic Sobolev space (see \emph{e.g.} \cite{troisi})
\[
W^{1,\overrightarrow{p}}(\Omega)=\left\{  u\in W^{1,1}(\Omega):
u_{x_i}\in L^{p_{i}}(\Omega),i=1,\ldots,N\right\}
\]
is a Banach space with respect to the norm
\begin{equation}
\left\Vert u\right\Vert _{W^{1,\overrightarrow{p}}(\Omega)}=\overset
{N}{\underset{i=1}{\sum}}\left\| u_{x_i}\right \| _{L^{p_{i}%
}(\Omega)}. \label{Sob_norm}%
\end{equation}
The space $W_{0}^{1,\overrightarrow{p}}(\Omega)$ is the closure of
$C_{0}^{\infty}(\Omega)$ with respect to the norm \eqref{Sob_norm}
and we will denote by $\left(
W_{0}^{1,\overrightarrow{p}}(\Omega)\right) ^{\prime}$ its dual.

\bigskip

%\subsection{Symmetrization}

A precise statement of our results requires the use of classical
notions of rearrangement and of suitable symmetrization of a Young
function, introduced by Klimov in \cite{Klimov 74}. \newline Let $u$
be a measurable function (continued by $0$ outside its domain)
fulfilling
\begin{equation}
\left\vert \{x\in\mathbb{R}^{N}:\left\vert u(x)\right\vert
>t\}\right\vert
<+\infty\text{ \ \ \ for every }t>0. \label{insieme livello di misura finita}%
\end{equation}
The \textit{symmetric decreasing rearrangement} of $u$ is the
function $u^{\bigstar}:\mathbb{R}^{N}\rightarrow\left[
0,+\infty\right[  $ $\ $satisfying
\begin{equation}
\{x\in\mathbb{R}^{N}:u^{\bigstar}(x)>t\}=\{x\in\mathbb{R}^{N}:\left\vert
u(x)\right\vert >t\}^{\bigstar}\text{ \ for }t>0. \label{livello palla}%
\end{equation}
The \textit{decreasing rearrangement} $u^{\ast}$ of $u$ is defined
by
\[
u^{\ast}(s)=\sup\{t>0:\mu_{u}(t)>s\}\text{ \ for }s\geq0,
\]
where
\[
\mu_{u}(t)=\left\vert \{x\in{\Omega}:\left\vert u(x)\right\vert
>t\}\right\vert \text{ \ \ \ \ for }t\geq0
\]
denotes the \textit{distribution function} of $u$. \newline
Moreover,
\[
u^{\bigstar}(x)=u^{\ast}(\omega_{N}\left\vert x\right\vert
^{N})\text{
\ \ }\hbox{\rm for $a.e.$}\;x\in{{\mathbb{R}}^{N}.}%
\]
Similarly, we define the \textit{symmetric increasing rearrangement}
$u_{\bigstar}$ on replacing \textquotedblleft$>$\textquotedblright\
by \textquotedblleft$<$\textquotedblright\ in the definitions of the
sets in (\ref{insieme livello di misura finita}) and (\ref{livello
palla}). We refer to \cite{BS} for details on these topics.

\bigskip

In this paper we will consider an $N-$\textit{dimensional Young
function} $\Phi:\rn \to \R$ (namely, an even convex function such
that $\Phi\left( 0\right)  =0$ \ \ and $\underset{\left\vert
\xi\right\vert \rightarrow+\infty}{\lim}\Phi\left(
\xi\right)  $ $=+\infty)$ of the following type:%
\begin{equation}
\Phi\left(  \xi\right)  =\underset{i=1}{\overset{N}{%
%TCIMACRO{\dsum }%
%BeginExpansion
{\displaystyle\sum}
%EndExpansion
}}\alpha_{i}\left\vert \xi_{i}\right\vert ^{p_{i}}\quad\text{ \ for
}\xi
\in\mathbb{R}^{N}\quad \text{ with }\alpha_{i}>0\quad \text{ for }i=1,...,N. \label{Fi pi}%
\end{equation}

We denote by $\Phi_{\blacklozenge}:\mathbb{R\rightarrow}\left[
0,+\infty \right[  $ the symmetrization of $\Phi$ introduced in
\cite{Klimov 74}. It is
the one-dimensional Young function fulfilling%
\begin{equation}
\Phi_{\blacklozenge}(\left\vert \xi\right\vert
)=\Phi_{\bullet\bigstar\bullet
}\left(  \xi\right)  \text{ \ for }\xi\in\mathbb{R}^{N}, \label{def fi rombo}%
\end{equation}
where $\Phi_{\bullet}$ is the Young conjugate function of $\Phi$ given by%
\[
\Phi_{\bullet}\left(  \xi^{\prime}\right)  =\sup\left\{
\xi\cdot\xi^{\prime
}-\Phi\left(  \xi\right)  :\xi\in\mathbb{R}^{N}\right\}  \text{ \ \ for \ }%
\xi^{\prime}\in\mathbb{R}^{N}.
\]
So $\Phi_{\blacklozenge}$ is the composition of Young conjugation,
symmetric increasing rearrangement and Young conjugate again.

We denote by $\overline{p}$  the \emph{harmonic average} of the
exponents $p_i$, \emph{i.e.}
\begin{equation}
\frac{1}{\overline{p}}=\frac{1}{N}\overset{N}{\underset{i=1}{\sum}}\frac{1}{p_{i}}\,.
\label{p barrato bis}%
\end{equation}
The \emph{harmonic average} $\overline{p}$ plays a basic role in
discussing anisotropic equations of the form \eqref{proto2}.
 Let us assume that $\overline{p}>1$ and set
\begin{equation}
\Lambda=\frac{2^{\overline{p}}\left(  \overline{p}-1\right)  ^{\overline{p}%
-1}}{\overline{p}^{\overline{p}}}\left[
\frac{ \displaystyle\prod_{i=1}^N p_{i}^{\frac{1}{p_{i}}}\left(  p_{i}^{\prime}\right)  ^{\frac{1}%
{p_{i}^{\prime}}}\Gamma(1+1/p_{i}^{\prime})}{\omega_{N}\Gamma(1+N/\overline
{p}^{\prime})}\right]  ^{\frac{\overline{p}}{N}}\left(
\prod_{i=1}^N\alpha_{i}^{\frac{1}{p_{i}}}\right) ^{\frac
{\overline{p}}{N}} \label{lapda}%
\end{equation}
with $\omega_{N}$ the measure of the $N-$dimensional unit ball,
$\Gamma$ the Gamma function and
$p_{i}^{\prime}=\frac{p_{i}}{p_{i}-1}$, the H\"{o}lder conjugate of
$p_i$ with the usual conventions if $p_{i}=1$. We are now in
position to evaluate $\Phi_{\blacklozenge}(\left\vert \xi\right\vert
)$. Easy calculations show (see \textit{e.g. }\cite{cianchi
anisotropo}) that
\begin{equation}
\Phi_{\blacklozenge}(\left\vert \xi\right\vert )=\Lambda\left\vert
\xi\right\vert ^{\overline{p}}. \label{fi rombo}%
\end{equation}
%where  the harmonic average of exponents $p_i$ defined in (\ref{p
%barrato bis}) fulfilling $\overline{p}>1$.
%and
%\begin{equation}
%\Lambda=\frac{2^{\overline{p}}\left(  \overline{p}-1\right)  ^{\overline{p}%
%-1}}{\overline{p}^{\overline{p}}}\left[
%\frac{ \displaystyle\prod_{i=1}^N p_{i}^{\frac{1}{p_{i}}}\left(  p_{i}^{\prime}\right)  ^{\frac{1}%
%{p_{i}^{\prime}}}\Gamma(1+1/p_{i}^{\prime})}{\omega_{N}\Gamma(1+N/\overline
%{p}^{\prime})}\right]  ^{\frac{\overline{p}}{N}}\left(
%\prod_{i=1}^N\alpha_{i}^{\frac{1}{p_{i}}}\right) ^{\frac
%{\overline{p}}{N}} \label{lapda}%
%\end{equation}
%with $\omega_{N}$ the measure of the $N-$dimensional unit ball,
%$\Gamma$ the Gamma function and
%$p_{i}^{\prime}=\frac{p_{i}}{p_{i}-1}$, the H\"{o}lder conjugate of
%$p_i$ with the usual conventions if $p_{i}=1$.
\newline In the anisotropic setting, we stress that $\overline{p}$ plays a
role also in a \textit{Polya-Szeg\"{o} principle} which reads as
follows (see \cite{cianchi anisotropo}). Let $u$ be a weakly
differentiable function in $ \rN$ satisfying (\ref{insieme livello
di misura finita}) and such that $$\sum_{i=1}^{N}\alpha_{i}
\int_{\rN} \left| u_{x_i} \right|^{p_{i}} dx<+\infty\,.$$ Then
$u^{\bigstar}$ is weakly differentiable in $\rN$ and
\begin{equation}
\Lambda\int_{\mathbb{R}^{N}}\left\vert \nabla
u^{\bigstar}\right\vert
^{\overline{p}}dx\leq\underset{i=1}{\overset{N}{%
%TCIMACRO{\dsum }%
%BeginExpansion
{\displaystyle\sum}
%EndExpansion
}}\alpha_{i}\int_{\mathbb{R}^{N}}\left\vert u_{x_i}
\right\vert ^{p_{i}}dx\text{ .\ } \label{polya}%
\end{equation}

\bigskip

\section{\textbf{Main results}}

%\subsection{Comparison result for elliptic problem}

In the present section, we focus our attention on the following
class of anisotropic elliptic problems
\begin{equation}
\left\{
\begin{array}
[c]{ll}%
-\operatorname{div}(a(x,u,\nabla u))+g(x,u)=f(x) & \hbox{in $\Omega$}\\
& \\
u=0 & \hbox{on $\partial\Omega$},
\end{array}
\right.  \label{problema_zero}%
\end{equation}
where $\Omega$ is a bounded open subset of $\R^{N}$ with Lipschitz
continuous boundary, $N\geq2$, $a:\Omega\times\mathbb{R}%
\times\mathbb{R}^{N}\rightarrow\mathbb{R}^{N}$ is a Carath\'{e}odory
function such that for \textit{a.e.} $x\in\Omega$, for all $s\in
\mathbb{R}$ and for all $ \xi, \xi' \in \rN$

\begin{itemize}
\item[(A1)] { $a(x,s,\xi)\cdot\xi\geq\overset{N}{\underset{i=1}{\sum}}%
\alpha_{i}\left\vert \xi_{i}\right\vert ^{p_{i}}\ \ \text{with }\alpha_{i}%
>0$ ,}

\item[(A2)] { $\left\vert a_{j}(x,s,\xi)\right\vert \leq\beta\left[
|s|^{\frac{\overline{p}}{p_{j}^{\prime}}}+|\xi_{j}|^{p_{j}-1}\right]
\text{ \ \ with }\beta>0$ \ \ $\forall j=1,\ldots,N$ ,}

\item[(A3)] {$\left(  a(x,s,\xi)-a(x,s,\xi^{\prime})\right)  \cdot\left(
\xi-\xi^{\prime}\right)  >0\qquad\hbox{for  $ \xi\neq \xi'$,}$}
\end{itemize}
 where $1\leq p_{1},\ldots,p_{N}<\infty$ are real numbers and
${{\overline{p}>1.}}$

\noindent
 Moreover, we assume that $g:\Omega\times
\R\rightarrow\mathbb{R}$ is a measurable, continuous and
non-decreasing function in $s$ for fixed $x$, and bounded in $x$
uniformly for bounded $u$ such that

\begin{itemize}
\item[(A4)]{ $g\left(  x,s\right)  s\geq b\left(  s\right)  s\qquad
\hbox{for \emph{a.e.}}\;\; x\in \Omega, \;\forall s\in \R$ , where
$b$ is a continuous and strictly increasing function such that
$b(0)=0.$}
\end{itemize}

\noindent
 Finally, we assume that
\begin{itemize}
\item[(A5)] $f:\Omega\rightarrow\mathbb{R}$ is a nonnegative function such
that {$f\in$}$\left(  W_{0}^{1,\overrightarrow{p}}(\Omega)\right)  ^{\prime}%
${.}
\end{itemize}

In order to give a precise statement of our results, we need to
precise what means to be less diffusive. Let $b_{1},b_{2}$ be two
{continuous strictly increasing functions. }We say that $b_{1}$ is
\textit{weaker} than
$b_{2}$ and we write%
\begin{equation}
b_{1}\prec b_{2}, \label{b minore b tilde}%
\end{equation}
if they have the same domains and there exists a
contraction\footnote{By contraction we mean $\left\vert
\rho(a)-\rho(b)\right\vert \leq\left\vert a-b\right\vert$ for $a, b
\in \R$.} $\rho: \mathbb{R}
\rightarrow%
%TCIMACRO{\U{211d} }%
%BeginExpansion
\mathbb{R}$ such that $b_{1}=\rho\circ b_{2}.$

We are interested in proving an integral estimate of a weak solution
$u\in W_{0}^{1,\overrightarrow{p}}(\Omega)$ to problem
\eqref{problema_zero} in terms of  the weak solution $w\in
W_{0}^{1,\overline{p}}(\Omega^{\bigstar})$ to the following problem
\begin{equation}
\left\{
\begin{array}
[c]{ll}%
-\operatorname{div}(\Lambda|\nabla w|^{\overline{p}-2}\nabla
w)+\widetilde
{b}(w)=\widetilde{f}(x) & \hbox{in $\Omega^\bigstar$}\\
& \\
w=0 & \hbox{on $\partial\Omega^\bigstar$},
\end{array}
\right.  \label{prob_sym_zero}%
\end{equation}
where $\Omega ^\bigstar$ is the ball centered at the origin and
having the same measure as $\Omega$,
\begin{itemize}
\item[(A6)] $\widetilde{b}$ {is a continuous and strictly increasing function such
that} $\widetilde{b}${$(0)=0$ ,}

\item[(A7)] $(\,\widetilde{b}\, )^{-1}\prec b^{-1}$,

\item[(A8)] $\widetilde{f}:\Omega^{\bigstar}\rightarrow\mathbb{R}$ is a
nonnegative radially symmetric function and decreasing along the
radii such that $\widetilde{f}\in\left(
W_{0}^{1,\overline{p}}(\Omega^{\bigstar })\right)  ^{\prime}.$
\end{itemize}

%We assume, in addition, that
%$$(\,\widetilde{b}\, )^{-1}\prec b^{-1}.$$

%In the
%next theorem we prove an $L^{\infty}-$norm estimate of the
%difference between the concentration of solutions to problems
%\eqref{problema_zero} and \eqref{prob_sym_zero} in terms of
%difference between the concentration of the corresponding data.
%Moreover, under the additional assumption that datum $f$ is less
%concentrated than datum $\widetilde{f},$ we obtain,in the corollary,
%a comparison result between the concentration of a solution $u$ to
%problem \eqref{problema_zero} and the concentration of the solution
%$w$ to problem (\ref{prob_sym_zero}).

\medskip

We stress that, by standard arguments and thanks to the results
contained in \cite{BB} (see also \cite{BCE} for the anisotropic
setting), there exists a unique weak
solution $w\in W_{0}^{1,\overline{p}}(\Omega^{\bigstar})$ to \eqref{prob_sym_zero} such that%
\begin{enumerate}
\item[$(i)$]{$\widetilde{b}(w)\in L^{1}\left(  \Omega^{\bigstar}\right)$}
\item[$(ii)$]{$\widetilde{b}(w)w\in L^{1}\left(  \Omega^{\bigstar}\right)$}
\item[$(iii)$]{$\Lambda\int_{\Omega}|\nabla w|^{\overline{p}-2}\nabla
w\cdot\nabla\phi \;dx+\int_{\Omega}\widetilde{b}(w)\,\phi\;dx\
=\langle\widetilde {f}\,,\phi\rangle_{\left(
W_{0}^{1,\overline{p}}(\Omega^{\bigstar })\right)  ^{\prime}}$}

for every $\phi\in W_{0}^{1,\overrightarrow{p}
}(\Omega^{\bigstar})\cap L^{\infty}(\Omega^{\bigstar})\text{ \ and }%
\varphi=w.$
\end{enumerate}

%
%\[
%\left\{
%\begin{array}
%[c]{l}%
%\widetilde{b}(w)\in L^{1}\left(  \Omega^{\bigstar}\right)  \text{ \
%and
%\ }\widetilde{b}(w)w\in L^{1}\left(  \Omega^{\bigstar}\right)  \text{ }\\
%\\
%\Lambda\int_{\Omega}|\nabla w|^{\overline{p}-2}\nabla
%w\cdot\nabla\phi \;dx+\int_{\Omega}\widetilde{b}(w)\,\phi\;dx\
%=\langle\widetilde {f}\,,\phi\rangle_{\left(
%W_{0}^{1,\overline{p}}(\Omega^{\bigstar
%})\right)  ^{\prime}}\;\quad\forall\phi\in W_{0}^{1,\overrightarrow{p}%
%}(\Omega^{\bigstar})\cap L^{\infty}(\Omega^{\bigstar})\text{ \ and }%
%\varphi=w.\text{\ }%
%\end{array}
%\right.
%\]
%

\bigskip

\begin{theorem}
\label{ellipticcomp} Assume that \emph{(A1)--(A8)} hold. Let $u$ be
a weak solution to the problem \eqref{problema_zero} and $w$  the
weak solution to the problem \eqref{prob_sym_zero}. Then,
\begin{equation}
\Vert(\mathcal{B}-\widetilde{\mathcal{B}})_{+}\Vert_{L^{\infty}(0,|\Omega
|)}\leq\Vert(\mathcal{F}-\widetilde{\mathcal{F}})_{+}\Vert_{L^{\infty
}(0,|\Omega|)}, \label{3:10}%
\end{equation}
where
\begin{align}
&  \mathcal{B}(s)=\int_{0}^{s}b(u^{\ast}(t))\;dt & \quad &
\widetilde
{\mathcal{B}}(s)=\int_{0}^{s}\widetilde{b}(w^{\ast}(t))\;dt\label{3:11}\\
&  \mathcal{F}(s)=\int_{0}^{s}f^{\ast}(t)\;dt & \quad &  \widetilde
{\mathcal{F}}(s)=\int_{0}^{s}\widetilde{f}^{\ast}(t)\;dt \label{3:11'}%
\end{align}
for $s\in(0,|\Omega|]$.
\end{theorem}

%%%%%%%%%%%%%%%%%%%%%%%%%%%%%%%%%%%%%%%%%%%%%%%%%%%%%%%%%%%%%%%%%%%%%%%%%%%%%%%%
%\begin{theorem}
%\label{ellipticcomp} Assume that (A1)- (A6) hold. If $w$ is a weak
%solution to problem \eqref{problema_zero}, then we have:
%\[
%\int_{0}^{s}w^{\ast}(\sigma)\,d\sigma\leq\int_{0}^{s}c^{\ast}(\sigma
%)\,d\sigma,\quad\forall s\in\lbrack0,|\Omega|],
%\]
%where $z$ is the solution to problem \eqref{prob_sym_zero}.
%\end{theorem}

\bigskip\ If we assume that the datum of problem \eqref{problema_zero}
dominates the datum of problem \eqref{prob_sym_zero}, then the
following comparison result between concentrations holds as an easy
consequence of Theorem \ref{ellipticcomp}.
\bigskip

\begin{corollary}
\label{C1}Under the same assumption of Theorem \ref{ellipticcomp},
if we suppose that
\begin{equation}
\mathcal{F}(s)\leq\widetilde{\mathcal{F}}(s)\qquad\hbox{for
any}\,\;s\in
\lbrack0,|\Omega|], \label{Dominazione f}%
\end{equation}
then
\begin{equation}
\mathcal{B}(s)\leq\widetilde{\mathcal{B}}(s)\qquad\hbox{for
any}\,\;s\in
\lbrack0,|\Omega|]. \label{3:12}%
\end{equation}
In particular, we have%
\begin{equation}\label{norm estim}
\int_{\Omega}\Psi(b(u(x)))\;dx\leq\int_{\Omega^{\bigstar}}\Psi(\widetilde
{b}(w(x)))\;dx \end{equation}
 for all convex and non-decreasing
function $\Psi:\mathbb{R}\rightarrow \mathbb{R}$ .
\end{corollary}

\bigskip

%\begin{remark}
\rm An immediate consequence of Corollary \ref{C1} are norm
estimates of $b(u)$ in terms of norm of $\widetilde{b}(w)$. An
example of applications of \eqref{norm estim} is the following one:
\[
\left\Vert b(u)\right\Vert _{L^{p}(\Omega)}\leq\left\Vert \widetilde
{b}(w)\right\Vert _{L^{p}(\Omega^{\bigstar})}\text{ \ \ for }1\leq
p\leq\infty\,.
\]
%and
%\[
%\left\Vert b(u)\right\Vert _{L^{p,q}(\Omega)}\leq\left\Vert
%\widetilde {b}(w)\right\Vert _{L^{p,q}(\Omega^{\bigstar})}\text{ \ \
%for }1\leq p<\infty\text{ and }1\leq q\leq\infty,
%\]
%where $L^{p,q}(\Omega)$ denotes the usual Lorentz space endowed with
%the norm
%\[
%\left\Vert v\right\Vert _{L^{p,q}(\Omega)}=\left\{
%\begin{array}
%[c]{ccc}%
%\left[
%%TCIMACRO{\dint _{0}^{\left\vert \Omega\right\vert }}%
%%BeginExpansion
%{\displaystyle\int_{0}^{\left\vert \Omega\right\vert }}
%%EndExpansion
%s^{\frac{1}{p}}\left(  \displaystyle\frac{1}{s}\int_{0}^{s}v^{\ast
%}(r)\;dr\right)^q  \;\displaystyle\frac{ds}{s}\right] ^{\frac{1}{q}}
%&
%\text{if} & 1<p<+\infty,1\leq q<\infty\\
%\underset{s\in\left(  0,\left\vert \Omega\right\vert \right)  }{\sup}%
%s^{\frac{1}{p}}\left(  \displaystyle\frac{1}{s}\int_{0}^{s}v^{\ast
%}(r)\;dr\right)  & \text{if} & 1<p<+\infty,q=\infty.
%\end{array}
%\right.
%\]
%
%
%\noindent Nevertheless, if $b$ is strictly increasing, then%
%\[
%\left\Vert u\right\Vert _{L^{\infty}(\Omega)}\leq\left\Vert
%w\right\Vert
%_{L^{\infty}(\Omega^{\bigstar})}\text{ .}%
%\]

%\end{remark}

\medskip

We emphasize that in the spirit of \cite{V}, Theorem
\ref{ellipticcomp} and Corollary \ref{C1} still hold if we do not
require the strictly monotony of $b$ and $\widetilde{b},$ but assume
that $b$ and $\widetilde{b}$ are {non-decreasing functions or, more
generally, maximal monotone graphs in $\r2$ such that $b(0)\ni0$ and
$\widetilde{b}\left(  0\right) \ni0$. Indeed, a maximal monotone
graph is a natural generalization of the concept of monotone
non-decreasing real function; moreover, the inverse of a maximal
monotone graph is again a maximal monotone graph (see \cite{V} for
more details).

\bigskip

%\begin{corollary}
%\label{C2}Under the same assumption of Theorem \ref{ellipticcomp}, if
%$-\operatorname{div}(a(x,u,\nabla u))=-\frac{\partial}{\partial x_{i}}\left(
%\left\vert \frac{\partial u}{\partial x_{i}}\right\vert ^{p_{i}-2}%
%\frac{\partial u}{\partial x_{i}}\right)  ,\widetilde{f}=$ $f^{\bigstar},$
%$\widetilde{b}=b,g=b,$ $b$ is a strictly increasing function and $\left\{
%x\in\Omega:w^{\bigstar}(x)>0\right\}  \subset\Omega^{\bigstar},$ then
%\[
%\left\vert \left\{  x\in\Omega^{\bigstar}:w^{\bigstar}(x)=0\right\}
%\right\vert \leq\left\vert \left\{  x\in\Omega:u(x)=0\right\}  \right\vert
%\]
%\end{corollary}
%\begin{corollary}
%\label{C3}Cambiare Under the same assumption of Theorem \ref{ellipticcomp}, if
%$\widetilde{f}=$ $f^{\bigstar},$ $\widetilde{b}=b,g=b$and b is a strictly
%increasing function and $\left\{  x\in\Omega:w^{\bigstar}(x)>0\right\}
%\subset\Omega^{\bigstar},$ then
%\[
%\left\vert \left\{  x\in\Omega^{\bigstar}:w^{\bigstar}(x)=0\right\}
%\right\vert \leq\left\vert \left\{  x\in\Omega:u(x)=0\right\}  \right\vert .
%\]
%\end{corollary}

%\section{Proofs}
%\subsection{Proof of Theorem \ref{ellipticcomp}}

\section{Proof of Theorem \ref{ellipticcomp}}

Let us consider the functions
$u_{\kappa,t}:\Omega\rightarrow$\textbf{ }$\mathbb{R}$ defined by
\[
u_{\kappa,t}\left(  x\right)  =\left\{
\begin{array}
[c]{ll}%
0 & \mbox{ if }\left\vert u\left(  x\right)  \right\vert \leq
t,\\
 &\\
\left(  \left\vert u\left(  x\right)  \right\vert -t\right)  \text{sign}%
\left(  u\left(  x\right)  \right)  & \mbox{ if
}t<\left\vert u\left(  x\right)  \right\vert \leq t+\kappa\\
& \\
\kappa\;\text{sign}\left(  u\left(  x\right)  \right)  & \mbox{ if
}t+\kappa<\left\vert u\left(  x\right)  \right\vert ,
\end{array}
\right.
\]
for any fixed$\ t$ and $\kappa>0.$ This function can be chosen  as a
test function in (\ref{problema_zero}). By (A1) and (A4),
\begin{equation}
-\frac{d}{dt}\int_{\{\left\vert u\right\vert >t\}}\overset{N}{\underset{i=1}{%
%TCIMACRO{\dsum }%
%BeginExpansion
{\displaystyle\sum}
%EndExpansion
}}\alpha_{i}\left\vert u_{x_i}\right\vert ^{p_{i}%
}dx\leq\int_{\{|u|>t\}}|f(x)|\;dx-\int_{\{|u|>t\}}b(u(x))\text{
sign}\;u\text{ }dx\qquad\hbox{for \emph{a.e.}
$t>0$}. \label{1111}%
\end{equation}
Taking into account \eqref{Fi pi}, \eqref{fi rombo} and
\eqref{polya}, analogous arguments as in \cite{cianchi anisotropo}
yield
\begin{equation}
-\frac{d}{dt}\int_{\{u^{\bigstar}>t\}}\Lambda|\nabla u^{\bigstar}|^{\overline{p}%
}dx\leq-\frac{d}{dt}\int_{\{\left\vert u\right\vert >t\}}\overset{N}%
{\underset{i=1}{%
%TCIMACRO{\dsum }%
%BeginExpansion
{\displaystyle\sum}
%EndExpansion
}}\alpha_{i}\left\vert u_{x_i}\right\vert ^{p_{i}%
}dx\qquad\hbox{for \emph{a.e.}
$t>0$}. \label{A}%
\end{equation}
%As consequence of (\ref{1111}),(\ref{A}) and Hardy-Littlewood
%inequality, we obtain
%\begin{equation}
%\label{E:40}0\leq-\frac{d}{dt}\int_{w^{\bigstar}>t}\Lambda|\nabla
%w^{\bigstar
%}|^{\overline{p}}\;dx\leq\int_{0}^{\mu_{w}(t)}f^{\ast}(s)\;ds-\int_{0}^{\mu_{w}%
%(t)}b\,(w^{\ast}(s))\;ds\qquad\hbox{for $t>0$}.
%\end{equation}
By the Coarea formula and the H\"{o}lder inequality,
\begin{equation}
\left(  -\frac{d}{dt}\int_{\{u^{\bigstar}>t\}}|\nabla u^{\bigstar}|^{\overline{p}%
}\;dx\right)  ^{\frac{1}{\overline{p}}}\geq
N\omega_{N}^{\frac{1}{N}}\mu _{u}(t)^{\frac{1}{N^{\prime}}}\left(
-\mu_{u}^{\prime}(t)\right) ^{-\frac
{1}{\overline{p}^{\prime}}}\qquad\hbox{for $ a.e.$ $t>0$} \label{E:41}.%
\end{equation}

\bigskip Since $f$ is nonnegative, the maximum principle
assures that $u\geq0$. Since $b$ is monotone, we obtain%
\begin{equation}
\int_{\{|u|>t\}}b(u(x))\text{ sign}\;u\text{
}dx=\int_{0}^{\mu_{u}(t)}b\,(u^{\ast
}(s))\;ds\qquad\hbox{for $ a.e. $ $t>0$}. \label{b}%
\end{equation}
Thus, as a consequence of \eqref{1111}, \eqref{A},\eqref{E:41} and
(\ref{b}), it follows that
\begin{equation}
\Lambda\left(  N\omega_{N}^{\frac{1}{N}}\mu_{u}(t)^{\frac{1}{N^{\prime}}%
}\left(  -\mu_{u}^{\prime}(t)\right)
^{-\frac{1}{\overline{p}^{\prime}}}\right)
^{\overline{p}}\leq\int_{0}^{\mu_{u}(t)}f^{\ast}(s)\;ds-\int_{0}^{\mu_{u}%
(t)}b\,(u^{\ast}(s))\;ds\qquad\hbox{for $ a.e. $ $t>0$}. \label{E:4}%
\end{equation}
%Putting
%\[
%\label{E:5}\mathcal{W}(s)=\int_{0}^{s}\lambda
%w^{\ast}(\sigma)\;d\sigma
%\quad\hbox{and}\quad\mathcal{F}(s)=\int_{0}^{s}f^{\ast}(\sigma)\;d\sigma
%\qquad\forall s\in\lbrack0,|\Omega|],
%\]
The relation \eqref{E:4} implies that
%\[
%\label{E:6}1\leq\frac{\left(  -\mu_{w}^{\prime}(t)\right)  ^{\frac{\overline{p}%
%}{\overline{p}^{\prime}}}\Lambda^{-1}}{\left(
%N\omega_{N}^{\frac{1}{N}}\mu _{w}(t)^{\frac{1}{N^{\prime}}}\right)
%^{\overline{p}}}\left[  \mathcal{F}(\mu
%_{w}(t))-\mathcal{B}(\mu_{w}(t))\right]  \qquad\hbox{for $t>0$},
%\]
%namely,%
\begin{equation}
1\leq\frac{-\mu_{u}^{\prime}(t)\Lambda^{-\frac{1}{\overline{p}-1}}}{\left(
N\omega_{N}^{\frac{1}{N}}\right)
^{\frac{\overline{p}}{\overline{p}-1}}\left( \mu
_{u}(t)\right)  ^{\frac{\overline{p}{\prime}}{N^{\prime}}}}\left[  \mathcal{F}%
(\mu_{u}(t))-\mathcal{B}(\mu_{u}(t))\right]  ^{\frac{1}{\overline{p}-1}}%
\qquad\hbox{for $ a.e. $ $t>0$}, \label{E:7}%
\end{equation}
where $\mathcal{F}$ and $\mathcal{B}$ are defined as in \eqref{3:11}
and \eqref{3:11'}, respectively.
%Now, integrating equation \eqref{E:7} between 0 and $t$, we have
%that
%\begin{equation}
%t\leq\left(  N\omega_{N}^{\frac{1}{N}}\right)  ^{-\overline{p}^{\prime}}%
%\Lambda^{-\frac{1}{\overline{p}-1}}\int_{\mu_{w}(t)}^{|\Omega|}\sigma^{-\frac
%{\overline{p}^{\prime}}{N^{\prime}}}\left[\mathcal{F}(\sigma)-
%\mathcal{B}(\sigma)\right]^{\frac{1}{\overline{p}-1}}\;d\sigma\qquad\hbox{for
%$t>0$},\label{E:8}%
%\end{equation}
%and so
%\begin{equation}
%w^{\ast}(s)\leq\left(  N\omega_{N}^{\frac{1}{N}}\right)  ^{-\overline{p}^{\prime}%
%}\Lambda^{-\frac{1}{\overline{p}-1}}\int_{s}^{|\Omega|}\sigma^{-\frac{\bar
%{p}^{\prime}}{N^{\prime}}}\left[\mathcal{F}(\sigma)-
%\mathcal{B}(\sigma)\right] ^{\frac{1}{\overline{p}-1}}\;d\sigma\qquad
%\hbox{for $s\in[0, |\Omega|]$}.\label{E:9}%
%\end{equation}
\newline By standard arguments (see, \emph{e.g.}, \cite{Ta1}), it follows that
\begin{equation}
\left(  -u^{\ast}(s)\right)  ^{\prime}\leq\left(  N\omega_{N}^{\frac{1}{N}%
}\right)
^{-\overline{p}^{\prime}}\Lambda^{-\frac{1}{\overline{p}-1}}s^{-\frac{\bar
{p}^{\prime}}{N^{\prime}}}\left[
\mathcal{F}(s)-\mathcal{B}(s)\right]
^{\frac{1}{\overline{p}-1}}\qquad\hbox{for \emph{a.e.} $s\in(0,
|\Omega|)$}. \label{E:10}%
\end{equation}
By \eqref{3:11},
\begin{equation}
\mathcal{B}^{\prime}(s)=b(u^{\ast}(s))\qquad\hbox{for
\emph{a.e.}}\;s\in(0,|\Omega|). \label{3:13}%
\end{equation}
Relations \eqref{3:11}, \eqref{E:10} and \eqref{3:13}  imply  that
\begin{equation}
\left\{
\begin{array}
[c]{lll}%
\Lambda\left(  N\omega_{N}^{\frac{1}{N}}\right)  ^{\overline{p}}s^{\frac{\overline{p}%
}{N^{\prime}}}\left[  -\displaystyle\frac{d}{ds}\left(  \gamma\left(
\mathcal{B}^{\prime}(s)\right)  \right)  \right]  ^{\overline{p}-1}+\mathcal{B}%
(s)\leq\mathcal{F}(s) & \qquad\hbox{for \emph{a.e.}
$s\in(0, |\Omega|])$} & \\
&  & \\
\mathcal{B}(0)=0,\;\;\mathcal{B}^{\prime}(|\Omega|)=0\,, &  &
\end{array}
\right.  \label{E:10'}%
\end{equation}
where $\gamma$ is the inverse function of $b$, \textit{i.e.} $\gamma
= b^{-1}.$

Let us consider problem \eqref{prob_sym_zero}. A weak solution $w$
to problem \eqref{prob_sym_zero} is unique and the symmetry of data
assures that $w(x)=w(|x|)$, \emph{i.e.} $w$ is positive and radially
symmetric. Moreover, setting $s=\omega_{N}|x|^{N}$ and
$\overset{\sim}{w}\left(  s\right) =w((s/\omega_{N})^{1/N})$, we get
that for all $s\in\lbrack0,|\Omega|]$
\[
-\Lambda|\overset{\sim}{w}^{\prime}\left(  s\right)  |^{\overline{p}-2}%
\overset{\sim}{w}^{\prime}\left(  s\right)  =\frac{s^{-\overline{p}/N^{\prime}}%
}{(N{\omega_{N}}^{1/N})^{\overline{p}}}\int_{0}^{s}\left(
f^{\ast}\left( \sigma\right)
-\widetilde{b}(\widetilde{w}(\sigma))\right)
\;d\sigma\qquad\hbox{for \emph{a.e.} $s\in(0, |\Omega|)$}.
\]
Since it is possible to show (see \cite[Lemma 1.31]{Dz85}) that the
above integral is positive, we deduce that $w(x)=w^{\bigstar}(x)$.
By the properties of $w$ we can repeat arguments used to prove
\eqref{E:10} replacing all the inequalities by equalities and
obtaining
\begin{equation}
\left(  -w^{\ast}(s)\right)  ^{\prime}=\left(  N\omega_{N}^{\frac{1}{N}%
}\right)
^{-\overline{p}^{\prime}}\Lambda^{-\frac{1}{\overline{p}-1}}s^{-\frac{\bar
{p}^{\prime}}{N^{\prime}}}\left[
\widetilde{\mathcal{F}}(s)-\widetilde {\mathcal{B}}(s)\right]
^{\frac{1}{\overline{p}-1}}\qquad\hbox{for \emph{a.e.} $s\in(0,
|\Omega|)$}\,. \label{E:10''}%
\end{equation}
Moreover, we have
\begin{equation}
\left\{
\begin{array}
[c]{lll}%
\Lambda\left(  N\omega_{N}^{\frac{1}{N}}\right)  ^{\overline{p}}s^{\frac{\overline{p}%
}{N^{\prime}}}\left[  \displaystyle-\frac{d}{ds}\left(
\widetilde{\gamma }\left( \widetilde{\mathcal{B}}^{\prime}(s)\right)
\right)  \right]
^{\overline{p}-1}+\widetilde{\mathcal{B}}(s)=\widetilde{\mathcal{F}}(s)
&
\qquad\hbox{for \emph{a.e.} $s\in(0, |\Omega|])$} & \\
&  & \\
\widetilde{\mathcal{B}}(0)=0,\;\;\widetilde{\mathcal{B}}^{\prime}%
(|\Omega|)=0, &  &
\end{array}
\right.  \label{PP}%
\end{equation}
where $\widetilde{\gamma}$ is the inverse function of
$\widetilde{b}$, \textit{i.e.} $\widetilde{\gamma} = (\,
\widetilde{b} \,)^{-1}.$

Since
$\mathcal{B},\widetilde{\mathcal{B}}\in\mathcal{C}([0,|\Omega|])$,
there exists $s_{0}\in(0,|\Omega|)$ such that
\begin{equation}
\Vert(\mathcal{B}-\widetilde{\mathcal{B}})_{+}\Vert_{L^{\infty}(0,|\Omega
|)}=(\mathcal{B}-\widetilde{\mathcal{B}})(s_{0}). \label{E:30}%
\end{equation}
In order to prove (\ref{E:10}), we argue by contradiction. Assume
that
\begin{equation}
(\mathcal{B}-\widetilde{\mathcal{B}})(s_{0})>\Vert(\mathcal{F}-\widetilde
{\mathcal{F}})_{+}\Vert_{L^{\infty}(0,|\Omega|)}\,. \label{E:31}%
\end{equation}
We distinguish  two cases: $s_{0}<|\Omega|$ and $s_{0}=|\Omega|.$

\textit{Case $s_{0}<|\Omega|$.} Combining (\ref{E:10'}) and
(\ref{PP}) yields
\begin{align}
&  \Lambda\left(  N\omega_{N}^{\frac{1}{N}}\right)
^{\overline{p}}s^{\frac{\bar {p}}{N^{\prime}}}\left[  \left(
-\frac{d}{ds}\left(  \gamma\left( \mathcal{B}^{\prime}(s)\right)
\right)  \right)  ^{\overline{p}-1}-\left(
-\frac{d}{ds}\left(  \widetilde{\gamma}\left(  \widetilde{\mathcal{B}}%
^{\prime}(s)\right)  \right)  \right)  ^{\overline{p}-1}\right] \label{E:35}\\
&  \leq\mathcal{F}(s)-\widetilde{\mathcal{F}}(s)+\widetilde{\mathcal{B}%
}(s)-\mathcal{B}(s) \qquad\hbox{for $ a.e. $ $s\in (0,
|\O|)$}\nonumber
\end{align}
By (\ref{E:31}),
\begin{equation}
\mathcal{F}(s)-\widetilde{\mathcal{F}}(s)+\widetilde{\mathcal{B}%
}(s)-\mathcal{B}(s)\leq\Vert(\mathcal{F}-\widetilde{\mathcal{F}})_{+}%
\Vert_{L^{\infty}(0,|\Omega|)}-(\mathcal{B}-\widetilde{\mathcal{B}})(s)<0
\label{E:36}%
\end{equation}
for $s\in(s_{0}-\varepsilon,s_{0}+\varepsilon).$ As a consequence of
(\ref{E:35}) and (\ref{E:36}) we obtain
\begin{align}
&  \Lambda\left(  N\omega_{N}^{\frac{1}{N}}\right)
^{\overline{p}}s^{\frac{\bar {p}}{N^{\prime}}}\left[  \left(
-\frac{d}{ds}\left(  \gamma\left( \mathcal{B}^{\prime}(s)\right)
\right)  \right)  ^{\overline{p}-1}-\left(
-\frac{d}{ds}\left(  \widetilde{\gamma}\left(  \widetilde{\mathcal{B}}%
^{\prime}(s)\right)  \right)  \right)  ^{\overline{p}-1}\right] \label{bbb}\\
&  =\Lambda\left(  N\omega_{N}^{\frac{1}{N}}\right)
^{\overline{p}}s^{\frac {\overline{p}}{N^{\prime}}}\omega(s)\left[
-\frac{d}{ds}\left(  \gamma\left( \mathcal{B}^{\prime}(s)\right)
-\widetilde{\gamma}\left(  \widetilde
{\mathcal{B}}^{\prime}(s)\right)  \right)  \right]  <0,\nonumber
\end{align}
where
\[
\omega(s)=\left(  \overline{p}-1\right)  \int_{0}^{1}\left\{  \left[
\tau\left( -\frac{d}{ds}\left( \gamma(\mathcal{B}^{\prime}(s)\right)
\right) +(1-\tau)\left( -\frac{d}{ds}\left( \widetilde{\gamma}\left(
\widetilde
{\mathcal{B}}^{\prime}(s)\right)  \right)  \right)  \right]  ^{\overline{p}%
-2}\right\}  \;d\tau>0.
\]
Setting
\begin{equation}
Z=\mathcal{B}-\widetilde{\mathcal{B}}\in
W^{2,\infty}(s_{0}-\varepsilon
,s_{0}+\varepsilon), \label{E:32}%
\end{equation}
we get
\begin{equation}
-\frac{d}{ds}\left(  \widetilde{\gamma}\left(
\mathcal{B}^{\prime}(s)\right) -\widetilde{\gamma}\left(
\widetilde{\mathcal{B}}^{\prime}(s)\right)
\right)  =-\frac{d}{ds}\left(  Z^{\prime}(s)\,\eta(s)\right)  , \label{aa}%
\end{equation}
where
\begin{equation}
\eta(s)=\int_{0}^{1}\widetilde{\gamma}^{\prime}\left(  \tau\mathcal{B}%
^{\prime}(s)+(1-\tau)\widetilde{\mathcal{B}}^{\prime}(s)\right)
\;d\tau>0.
\label{E:34}%
\end{equation}
By (A7), we can conclude that

\begin{equation}
-\frac{d}{ds}\left(  \gamma\left(  \mathcal{B}^{\prime}(s)\right)
-\widetilde{\gamma}\left(  \mathcal{B}^{\prime}(s)\right)  \right)
\geq0 \qquad\hbox{for $ a.e. $ $s\in (0, |\O|)$}.
\label{ccc}%
\end{equation}
Then, by (\ref{aa}) and \eqref{ccc},
\begin{equation}
-\frac{d}{ds}\left(  Z^{\prime}(s)\,\eta(s)\right)
\leq-\frac{d}{ds}\left( \gamma\left(  \mathcal{B}^{\prime}(s)\right)
-\widetilde{\gamma}\left(
\widetilde{\mathcal{B}}^{\prime}(s)\right)  \right) \qquad\hbox{for $ a.e. $ $s\in (0, |\O|)$} . \label{E:33}%
\end{equation}
Finally, thanks to (\ref{bbb}) and (\ref{E:33}), we have
\begin{align}
&  \Lambda\left(  N\omega_{N}^{\frac{1}{N}}\right)
^{\overline{p}}s^{\frac{\bar {p}}{N^{\prime}}}\omega(s)\left(
-\frac{d}{ds}\left(  \eta(s)Z^{\prime
}(s)\right)  \right)  \leq\label{E:37}\\
&  \leq\Lambda\left(  N\omega_{N}^{\frac{1}{N}}\right)  ^{\overline{p}}%
s^{\frac{\overline{p}}{N^{\prime}}}\omega(s)\left[
-\frac{d}{ds}\left( \gamma\left(  \mathcal{B}^{\prime}(s)\right)
-\widetilde{\gamma}\left( \widetilde{\mathcal{B}}^{\prime}(s)\right)
\right)  \right]  <0 \qquad\hbox{for $ a.e. $ $s\in (0,
|\O|)$}.\nonumber
\end{align}
We can conclude that
\begin{equation}
-\frac{d}{ds}\left(  \eta(s)Z^{\prime}(s)\right)  <0\qquad\hbox{for
}\;s\in(s_{0}-\varepsilon,s_{0}+\varepsilon), \label{ZZ}%
\end{equation}
which is in contradiction with the assumption (\ref{E:30}),
\textit{i.e.} $Z$ has a maximum in $s_{0}.$

\textit{Case $s_{0}=\left\vert \Omega\right\vert $.} In this case,
the inequality (\ref{ZZ})\ holds for $s\in(\left\vert
\Omega\right\vert -\varepsilon,\left\vert \Omega \right\vert ).$ So
$Z^{\prime}(\left\vert \Omega\right\vert )>0,$ but this is not true
since $Z^{\prime}(\left\vert \Omega\right\vert )=0.$

\qed

\bigskip

%\subsection{\bigskip Proof of Corollaries}
%Proof of Corollary \ref{C2}
%\bigskip$\left\{  x\in\Omega:w^{\bigstar}(x)>0\right\}  \subset\Omega
%^{\bigstar}$
%If $\widetilde{f}=$ $f^{\bigstar},$ $\widetilde{b}=b,$ $b$ is a strictly
%increasing function by Theorem \ref{ellipticcomp} it follows that
%\[
%\int_{0}^{s}b(u^{\ast}(t))\;dt\leq\int_{0}^{s}\widetilde{b}(w^{\ast
%}(t))\;dt\text{ \ for every }s\in\left[  0,|\Omega|\right]
%\]
%Moreover
%\begin{align*}
%\int_{0}^{\left\vert \Omega\right\vert }b(u^{\ast}(t))\;dt &  =\int_{\Omega
%}b(u(x))\;dx=\int_{\Omega}\left[  f(x)+\frac{\partial}{\partial x_{i}}\left(
%\left\vert \frac{\partial u}{\partial x_{i}}\right\vert ^{p_{i}-2}%
%\frac{\partial u}{\partial x_{i}}\right)  \right]  dx\\
%&  =
%\end{align*}%
%\[
%\int_{0}^{\Omega}b(u^{\ast}(t))\;dt\leq\int_{0}^{\Omega^{\bigstar}}%
%\widetilde{b}(w^{\ast}(t))\;dt.
%\]
%It follows that
%\[
%\int_{0}^{s}b(u^{\ast}(t))\;dt=\int_{0}^{\Omega}b(u^{\ast}(t))\;dt\geq\int
%_{0}^{\Omega^{\bigstar}}b(w^{\ast}(t))\;dt
%\]
%and $\left\{  x\in\Omega:w^{\bigstar}(x)>0\right\}  \subset\Omega^{\bigstar},$
%then
%\[
%\left\vert \left\{  x\in\Omega^{\bigstar}:w^{\bigstar}(x)=0\right\}
%\right\vert \leq\left\vert \left\{  x\in\Omega:u(x)=0\right\}  \right\vert
%\]

\paragraph*{Acknowledgements}

This work has been partially supported by GNAMPA of the Italian
INdAM (National Institute of High Mathematics) and ``Programma
triennale della Ricerca dell'Universit\`{a} degli Studi di Napoli
``Parthenope'' - Sostegno alla ricerca individuale 2015-2017''.

\end{document}